\documentclass{amsart}

\usepackage{amsmath, amsfonts, amssymb, amscd}

\newtheorem{theo}{Theorem}[section]
\newtheorem{lem}[theo]{Lemma}
\newtheorem{prop}[theo]{Proposition}

\newtheorem{cor}[theo]{Corollary}
\newtheorem{rem}[theo]{Remark}

\newtheorem{definition}[theo]{Definition}
\newenvironment{pf}{\noindent{\it Proof.
}}{$\blacksquare$\par\medskip}
\newenvironment{pfnqed}{\noindent{\it Proof.
}}{\par\medskip}

\newcommand{\C}{{\mathbb C}}
\newcommand{\R}{{\mathbb R}}

\newcommand{\g}{{\mathfrak g}}

\renewcommand{\k}{{\mathfrak k}}
\newcommand{\m}{{\mathfrak m}}

\newcommand{\z}{{\mathfrak z}}

\newcommand{\h}{{\mathfrak h}}
\newcommand{\p}{{\mathfrak p}}

\renewcommand{\t}{{\mathfrak t}}

\newcommand{\B}{{\mathcal B}}

\renewcommand{\L}{{\mathcal L}}
\newcommand{\CC}{{\mathcal C}}
\newcommand{\Fut}{{\mathcal F}}
\newcommand{\CU}{{\mathcal{U}}}

\newcommand{\ch}[1]{{#1}^{\vee}}
\renewcommand{\varpi}{\psi}

\newcommand{\co}{m}

\newcommand{\Aut}{\operatorname{Aut}}
\newcommand{\aut}{\mathfrak{aut}}
\newcommand{\Ad}{\operatorname{Ad}}

\newcommand{\reg}{{\operatorname{reg}}}

\renewcommand{\=}{\overset{\text{def}}{=}}

\def\sideremark#1{\ifvmode\leavevmode\fi\vadjust{% The remark
\vbox to0pt{\hbox to 0pt{\hskip\hsize\hskip1em% will appear only
\vbox{\hsize3cm\tiny\raggedright\pretolerance10000% on the side
\noindent #1\hfill}\hss}\vbox to8pt{\vfil}\vss}}}% in 3cm

\title[K\"ahler-Ricci solitons on homogeneous toric bundles  (II)]
{K\"ahler-Ricci solitons \\ on homogeneous toric bundles  (II)}
\author{Fabio Podest\`a
and Andrea Spiro}


\begin{document}

  %\vspace*{.1cm}
\begin{abstract}
It is proved that an homogeneous toric bundles over
a flag manifold $G^\C/P$ admits a  K\"ahler-Ricci solitonic metric if and only if it is Fano.
In particular,  an  homogeneous toric bundle of this kind is K\"ahler-Einstein
if and only if it is Fano and its Futaki invariant  vanishes identically.
\end{abstract}

\maketitle

%\null
% \vspace*{-.3cm}

\section{Introduction}
\bigskip
In this paper we continue the discussion of  \cite{PS2} on K\"ahler-Einstein and  K\"ahler-Ricci solitonic metrics
over  homogeneous bundles $\pi: M \to V$, with fiber equal to   a compact toric K\"ahler manifold $F$
and basis $V$ equal to   a
generalized flag manifold $V = G^\C/P$ of a complex
semisimple Lie group $G^\C$.  We
  call any such bundle a {\it homogeneous toric
bundle}. \par
  \smallskip
In \cite{PS2} we  gave necessary and sufficient conditions in order that a homogeneous toric bundle $\pi: M \to V = G^\C/P$ has positive first Chern class; In particular this occurs only if $F$ is Fano. In this second
part we determine when a homogenous toric bundle admits a  K\"ahler-Ricci soliton.\par
We recall that a {\it K\"ahler-Ricci soliton\/}ÿ consists of  a K\"ahler form $\omega$ associated with a  (real) vector field $X$ such that
$$\rho - \omega = \L_X \omega\ ,\qquad \L_{J X} \omega = 0\ ,$$
where $\rho$ denotes the Ricci form of $\omega$. Notice that  if the associated vector field $X$ is trivial,
the  K\"ahler-Ricci soliton $\omega$ is a K\"ahler-Einstein form.\par
\medskip
Our main result is the following.\par
\medskip
\begin{theo} \label{maintheorem}   Let $F$ be a toric K\"ahler manifold of dimension $m$ and
 $\pi: M \to V$ be a homogeneous toric bundle with fiber $F$ and basis $V = G^\C/P$. The bundle
 $ M$ admits a K\"ahler-Ricci soliton  if and only if it is Fano.
\par
 In particular,  the bundle $M$ is K\"ahler-Einstein if and only if it is Fano and
 its Futaki functional vanishes identically.
\end{theo}
\bigskip
This theorem  extends the result of X.-J. Wang and X. Zhu (\cite{WZ}) who proved the existence of a K\"ahler-Ricci soliton on any Fano toric manifold $F$, i.e. when the basis of the toric bundle reduces to
a single point. On the other hand, our theorem includes the results of  N. Koiso and Y. Sakane  in \cite{Sa, KS, KS1}, which give necessary and sufficient conditions for
homogeneous toric bundles with fiber $\C P^1$  in order to be  K\"ahler-Einstein (see also \cite{PS, DW}).  It also generalizes  Koiso's result  (\cite{Ko}) on the existence of a K\"ahler-Ricci soliton on  any Fano, homogeneous toric bundle  with fiber $\C P^1$ (see also \cite{TZ1}).\par
\medskip
The paper is organized as follows. In \S 2 we fix notations and recall some facts on homogeneous toric bundles that
were used and/or proved in \cite{PS2}. In \S 3 we recall  the fundamental results of \cite{TZ,TZ1} on K\"ahler-Ricci
solitons and we obtain some consequences on homogeneous toric bundles. In \S 4 we compute the holomorphic invariant
introduced by Tian and Zhu in \cite{TZ1} and in \S 5 we show that the problem of finding a K\"ahler-Ricci soliton on
the homogeneous toric bundles can be reduced to a suitable partial differential equation on the toric fiber $F$:
This equation turns out to be  very close to the equation studied in \cite{WZ}. We conclude showing that under
suitable modifications, the arguments used in the proof of Wang and Zhu  for the solvability of that equation works
in our case as well.\par
We remark that from the proof of Theorem \ref{maintheorem}, it follows that a vector field $X$ on a homogeneous toric
bundle over a flag manifold is the associated vector field of a K\"ahler-Ricci soliton  if and only if the Tian and
Zhu's invariant $\Fut_X(\cdot)$ vanishes identically.

\bigskip
\section{Notations and preliminaries}
 \bigskip
 As we mentioned in the Introduction, this paper is the continuation of  \cite{PS2}
and we will constantly use  the same notation and definitions introduced  in that paper. For readers convenience, we  briefly recall here all  notations and definitions adopted in that paper, but we refer to \cite{PS2}
for more detailed information.\par
\bigskip
 For any Lie group $G$, we will
denote its Lie algebra by the corresponding gothic letter $\g$.
Given a Lie homomorphism $\tau:G\to G'$, we will always use the
same letter to represent the induced Lie algebra homomorphism
$\tau:\g\to\g'$.  The center of $G$ will be denoted by $Z(G)$ and the center of $\g$ by $\z(\g)$.\par
If $G$ acts on a manifold $N$, for any $X\in
\g$, we will use the symbol $\hat X$ to indicate
  the corresponding induced vector field on $N$. We recall here that
$\widehat{[X,Y]} = - [\hat X,\hat Y]$ for every
  $X,Y\in \g$.\par
 We will also denote by $N_\reg$ the set of
  $G$-principal points in $N$.\par
\medskip
  The Cartan Killing form of a semisimple Lie algebra $\g$
will be always denoted by $\B$ and, for any
  $X\in \g$, we set $\ch X= - \B(X, \cdot) \in \g^*$.
Given a root system $R$ w.r.t. a
fixed maximal torus, we will denote by $E_\alpha\in\g^\C$
the root vector corresponding to the root $\alpha$ in
the Chevalley normalization and by
$H_\alpha=[E_\alpha,E_{-\alpha}]$ the $\B$-dual of
$\alpha$.\par
\bigskip
In all the following, $F$ denotes  a compact, toric K\"ahler
manifold  with $\dim_\C F = m$ and we indicate by $T^m$ the $m$-dimensional torus
acting effectively on $F$ by holomorphic isometries. A homogeneous   toric bundle
 is a compact K\"ahler manifold of the form
$$M = G^\C \times_{P, \tau}  F = G \times_{K, \tau}  F\eqno(2.1)$$
where $V = G^\C/P = G/K$ is a flag manifold of (complex) dimension $n$, $G$ is a compact semisimple
Lie group, $G^\C$ its complexification, $P$ a suitable parabolic subgroup  and $\tau: P \to (T^m)^\C$ is a surjective homomorphism. \par
\smallskip
We will  constantly identify $F$ with the fiber $F = F_{e K} = \pi^{-1}(eK)$ over the base point
$eK \in V = G/K$.\par
\smallskip
The complex structures of $M$, $F$ and $V$ will be denoted by $J$, $J_F$ and $J_V$, respectively.
Notice that $J_V$ is the natural $G^\C$-invariant complex structure of the complex homogeneous space $G^\C/P$ and that $J$ is the unique $G^\C$-invariant complex structure on $M$, which makes $\pi: M \to V$ a holomorphic map and induces on $F = \pi^{-1} (e K)$ the complex structures $J_F$.\par
\par
\smallskip
 We observe that
both $G^\C$ and  $(T^\co)^\C$ act naturally  as groups of holomorphic transformations on $(M, J)$, with
two commuting actions. The action of $G^\C$ is the one induced on $M$  by its  standard action on $G^\C \times   F$,
while the action of $(T^\co)^\C$   is defined by
$$h([g, x]_{K,\tau}) \= [ g, h^{-1}(x)]_{K,\tau}\ ,\qquad \text{for any}\ h\in (T^m)^\C\ .$$
For this reason, in the following  we will identify  $G^\C \times (T^m)^\C$ with the corresponding
subgroup of  $\Aut(M, J)$ and $\g^\C + \t^\C$ will be identified with the corresponding  subalgebra of $\aut(M,J) = Lie(\Aut(M,J))$.\par
\medskip
We recall that   $\g$  admits an
$\Ad(K)$-invariant
  decomposition $\g=\k \oplus \m$ and that,
for any fixed CSA $\h \subset \k^\C$ of $\g^\C$, the
associated root system $R$ admits a corresponding decomposition
$R = R_o + R_\m$, so that $E_\alpha \in \k^\C$ if $\alpha \in R_o$
and $E_\alpha \in\m^\C$ if $\alpha\in R_\m$.
Furthermore, $J_V$ induces a splitting $R_\m = R_\m^+\cup R_\m^-$
into two disjoint subset of positive and negative roots, so that
the $J_V$-holomorphic and $J_V$-antiholomorphic subspaces of
$\m^\C$ are given by
$$\m^{(1,0)} = \sum_{\alpha\in R_\m^+}\C E_\alpha,\quad
\m^{(0,1)} = \sum_{\alpha\in R_\m^-}\C
E_\alpha\ .\eqno(2.2)$$
The Lie algebra $\p$ of the parabolic subgroup $P$ is $\p
= \k^\C + \m^{(0,1)}$. \par
 We also recall that for any $G$-invariant K\"ahler
form $\omega$ of $V$ there
exists a uniquely associated
  element $Z_\omega \in \z(\k)$ so that
$\left.\omega(\hat X, \hat Y)\right|_{eK} = \B(Z_\omega,
[X,Y])$
for any $X, Y\in \g$. In particular, the
$G$-invariant K\"ahler-Einstein form $\omega_V$ on $V$,
with Einstein constant $c =1$, is associated with the element
$$Z_V = - \frac{1}{2 \pi}ÿ\sum_{\alpha \in R^+_\m}
i H_{\alpha}\ .\eqno(2.3)$$
(see e.g. \cite{Be, BFR} - be aware that in this paper,
we adopt the definition of Ricci form $\rho$ used e.g. in \cite{Fu1}, which differs from
the one in  \cite{Be} and
\cite{BFR}   by the factor $\frac{1}{2 \pi}$).\par
\medskip
The homomorphism $\tau: P\to (T^\co)^\C$ is completely determined by its restriction to the connected component of the identity $Z^o(K)$ of $Z(K)$, which gives   a surjective homomorphism $\tau: Z^o(K) \to T^m$ and a surjective Lie algebra homomorphism $\tau: \z(\k) \to \t$.  \par
In the following, we will denote by  $\t_G \= (\ker \tau)^\perp \cap \z(\k)$. Notice that $\t_G$ integrates to a closed subtorus and moreover we can choose a $\B$-orthonormal basis
$(Z_1,\ldots,Z_\co)$ of $\t_G$ so that $\exp(\R\cdot Z_j)$ is closed for every $j=1,\ldots,\co$. We will denote by  $Z_j' \= \tau(Z_j)$ for $j=1,\ldots,\co$ and by $\nu_j$ the smallest real number such that
$\exp(\nu_j Z'_j) = e$.
\par
We will also denote by $(F_{\alpha_1}, G_{\alpha_1}, \dots, F_{\alpha_n}, G_{\alpha_n})$ the basis
for  $\m  \subset \g$ given by the element
$$F_{\alpha_i} = \frac{1}{\sqrt{2}} (E_{\alpha_i} - E_{- \alpha_i})\ ,\qquad
G_{\alpha_i} = \frac{i}{\sqrt{2}} (E_{\alpha_i} + E_{- \alpha_i})\ ,\qquad \alpha_i \in R_\m\ .\eqno(2.4)$$
\par
\bigskip
The rest of this section will be devoted to  the properties  we will later use of the so-called  ``algebraic representative" of
a closed $G$-invariant  2-form and of their relations with the moment maps. \par
\smallskip
If $\varpi$ is a $G$-invariant closed 2-form on $M$, then
there exists a  $G$-equivariant map $Z_\varpi: M \to \g$, uniquely associated
with $\varpi$ so that
$$\varpi_p(\hat X, \hat Y) = \B([Z_\varpi|_p, X], Y) =
\B(Z_\varpi|_p, [X, Y])\qquad \text{for any}\ X, Y \in
\g\eqno(2.5)$$
This map is called   {\it algebraic representative of\/} $\varpi$ and, in case
$\varpi$ is non-degenerate, the moment map  determined by
$\varpi$
$$\mu_\varpi: M \to \g^*$$
coincides with  the $(-\B)$-dual map  of $Z_\psi$
$$\ch Z_\varpi \= - \B(Z_\varpi, \cdot): M \to \g^*.$$
\par
\medskip
By $G$-equivariance,  any algebraic
representative $Z_\varpi$ is uniquely determined
by its restriction on  the fiber $F= \pi^{-1}(eK)$ and such restriction $Z_\varpi|_F$ takes values  in  $\z(\k)$. \par
In case a $G$-equivariant 2-form  $\varpi$ is cohomologous to $0$,  its restriction to $F$  must be of the form
$\varpi = d d^c \phi$ for some $K$-invariant smooth function $\phi: F \to \R$ and the restriction to $F$ of its
algebraic representative is
$$Z_\psi|_F = - \sum_i J \hat Z_i(\phi)Z_i\ .\eqno(2.6)$$
For any given K\"ahler form $\omega \in c_1(M)$, the restrictions to $F$ of the algebraic representatives of
$\omega$ and of its Ricci form
$\rho$ are  as follows:
$$Z_\omega|_F = \sum_i f_i Z_i + Z_V\ ,\qquad \text{for some smooth functions}\ \ f_i: F \to \R\ ,\eqno(2.7)$$
$$Z_\rho|_F =  \sum_{i = 1}^\co \frac{J \hat Z_i(\log h)}{4 \pi}
Z_i +
Z_V\ ,\eqno(2.8)$$
where $Z_V \in \z(\k)$ is
the element defined in (2.3)  and
$$h =  \det\left(\begin{matrix}
- f_{i,j}\end{matrix}\right)\cdot \prod_{\alpha \in R^+_\m}
( a^i_{\alpha} f_i + b_\alpha)\ ,\eqno(2.9)$$
where $f_{i,j} \= J \hat Z_j(f_i)$,
$a^i_{\alpha} \= \alpha(i Z_i)$,
  $b_\alpha \= \alpha(i Z_V)$. \par
Moreover, for any $p \in F$,
$$- f_{i,j}(p) =  \omega_p(\hat Z_j, J \hat Z_i) = \frac{1}{2 \pi} g_p(\hat Z_i, \hat Z_j)  \ .\eqno(2.10)$$
In other words, for any point $p \in F$, the values $-f_{i,j}(p)$ are the entries  of a symmetric, positive
definite matrix  and  one can check that  the map
$$\mu: F \to \t^*\ ,\qquad \mu(q) = -\left.\B(\sum_\ell f_\ell Z_\ell, \cdot)\right|_{\t_G}\in \t_G^* \simeq \t
\eqno(2.11)$$
is a moment map for the action of $T^m$  determined by $\omega|_{TF}$. \par

\smallskip
If $F$ has $c_1(F) > 0$, for a given $T^m$-invariant
 K\"ahler form $\psi \in c_1(F)$,
a corresponding moment map $\mu_\psi: F \to \t^*$ is called {\it metrically normalized\/}  if $\int_F \psi \cdot  \eta_\psi^\co = 0$ where $\eta_\psi$ is the unique K\"ahler form in $c_1(F)$ that has $\psi$ as Ricci form.  By \cite{PS2} \S 4,  for any $\psi \in c_1(F)$, there exists a unique associated  metrically normalized moment map and the polytope
$\Delta_F = \mu_{\psi}(F)$ is independent of $\psi$ and it is called the {\it canonical polytope of $F$\/}.\par
If $M$ is Fano, by Thm. 1.1 of \cite{PS2}, then also  $F$ is Fano  and the moment
map defined in (2.11) is the metrically normalized moment map determined by $\omega|_{TF}$. In particular, $\mu(F) = \Delta_F$.\par
 \medskip
In all the following,   we will also constantly  identify $\t^* (= \t^*_G)$ with $\R^\co$, through the vector space isomorphism that
maps  the elements $ e_\ell \= \frac{1}{4 \pi} \left.\B(Z_\ell, \cdot)\right|_{\t_G}$ into  the canonical basis of $\R^\co$.
By virtue of  such identification,  in next sections the map $\mu$ will always  be written as
$$\mu: F \to \Delta_F \subset \R^\co\ ,\qquad \mu(q) =  \left (- 4 \pi f_1(q), \dots, - 4 \pi f_\co(q)\right)\ .\eqno(2.12)$$
\bigskip
\bigskip
\section{K\"ahler-Ricci solitons and associated vector fields on  homogeneous toric bundles}
\bigskip
First of all, let us recall the definition of K\"ahler-Ricci soliton (see e.g. \cite{TZ}).\par
\medskip
\begin{definition} \label{solitondef}  {\rm Let $(N, J, \hat \omega)$
 be a compact K\"ahler manifold of positive first Chern class. We call {\it K\"ahler-Ricci soliton\/} any pair $(\omega, X)$, where $\omega$ is a K\"ahler form on $M$ and $X$ is a (real) vector field  on $M$ such that:
$$a)\ \L_{JX} \omega = 0\ ,\qquad b)\  \rho - \omega = \L_X \omega = d(\imath_X \omega)\ .$$
If $(\omega, X)$ is a K\"ahler-Ricci soliton, we will say that  $\omega$ is the  {\it K\"ahler form of the soliton\/}  and  that
$X$ is the   {\it  associated vector field\/}.}
\end{definition}
\bigskip
From b) it is clear that a compact K\"ahler manifold admits a K\"ahler-Ricci soliton only if it
is Fano.\par
\medskip
We need now to introduce some notation regarding K\"ahler-Ricci solitons. For a given compact K\"ahler manifold $(N,J)$, we
 denote by $\Aut(N, J)$ the
group of all complex automorphisms of $N$, by $\Aut(N, J)^o$ its connected component of the identity and by
$R_u(N,J)$ its unipotent radical. \par
\medskip
In the next statement, we collect some crucial  facts on K\"ahler-Ricci solitons obtained  by Tian and Zhu (see \cite{TZ}, Thm. A,  \cite{TZ1}, Prop. 3.1, Prop. 2.1, Thm. 3.2).
\medskip
\begin{theo} \label{TZtheorem} Let $(N, J, \hat \omega)$ be a compact K\"ahler manifold with positive first Chern class and assume that it admits
a K\"ahler-Ricci soliton $(\omega, X)$.  Denote also by $G^{(\omega)}\subset \Aut(N,J)^o$ the
subgroup of all isometries of $\omega$. Then:
\begin{itemize}
\item[i)] $G^{(\omega)}$ is a maximal compact subgroup of $\Aut(N,J)^o$ and $JX$ belongs to the center $\z(\g^{(\omega)})$;
\item[ii)]  all K\"ahler-Ricci solitons $(\omega', X')$ on $N$ are of the form
$$\omega' = \sigma^* \hat \omega\ ,\qquad X' = \sigma^{-1}_*(X)$$
for some $\sigma \in \Aut(N,J)^o$.
\end{itemize}
\end{theo}
\bigskip
\begin{rem} \label{remark2.3} {\rm From i) and ii) of the previous theorem, it follows immediately that
{\it $N$ admits a K\"ahler-Ricci soliton $(\omega, X)$ if and only if,  for a given maximal compact subgroup $\tilde  G \subset \Aut(N,J)^o$, there is   a  K\"ahler-Ricci soliton $(\omega^{(\tilde G)}, X^{(\tilde G)})$,  where $\omega^{(\tilde G)}$ is $\tilde G$-invariant and
  $J X^{(\tilde G)}   \in  \z(\tilde \g)$\/}}.
\end{rem}
\bigskip
Let us now consider an  homogeneous toric bundle $M= G^\C \times_{P, \tau}  F$, with the fiber $F$ acted on by the torus $T^\co$ (and hence by its complexification $(T^\co)^\C$).
The following lemma  is crucial.\par
\medskip
\begin{lem} \label{center} Let $M= G^\C \times_{P, \tau}  F$ be a homogeneous
toric bundle, $\tilde G \subset \Aut(M)^o$ a maximal compact subgroup
containing $G\times T^\co$ and $\tilde \g = Lie(\tilde G)$.
Then  $\z(\tilde \g) \subset \t $.
\end{lem}
\begin{pf}
By  Blanchard's Lemma (\cite{Bl}; see also \cite{Ah}, Prop. 1, p. 45), for any $Y \in \z(\tilde \g)$ the group
$A \= \overline{\exp(\R \cdot Y)}$ is a  compact, abelian  subgroup of  $Z(\tilde G)$ consisting of  fiber preserving biholomorphisms.
This implies that $A$
 projects onto  a  compact, connected group $A_V$ of biholomorphisms
of $V$  with   $A_V \subseteq C_{\mathcal A}(G^\C)$, where ${\mathcal A} \= \Aut(V, J_V)^o$. Now,
if $\mathcal A = G^\C$, then $A_V$ is trivial because $G^\C$ is semisimple.  If
 $\mathcal A \supsetneq G^\C$, then the possible pairs $({\mathcal A}, G^\C)$ have been classified by Onishchik in \cite{On}
 and it is easily checked that $C_{\mathcal A}(G^\C)$ is trivial and hence $A_V = \{e\}$. This means that $A$ fixes all fibers
 and that  the restriction of $A$ to
 $F = F|_{eK}$ commutes with the action of $T^\co$. \par
 On the other hand, by Demazure's Structure Theorem for  toric manifold (see e.g. \cite{Od}, p. 140),  $\Aut(F, J_F)^o$ is a linear algebraic group and $T^\co$
  is a maximal algebraic torus of $\Aut(F, J_F)^o$. This implies that for any $a \in A$ the biholomorphism $a|_{F } : F \to F$  coincides with some  biholomorphism
  $t|_F$, $t \in T^m$ or, equivalently, that $a \circ t^{-1}|_{F} = Id$.
Since both $a$ and $t$  commute with  $G$, it  follows that $a \circ t^{-1}|_{\pi^{-1}(gK)} = Id$ for any
fiber $\pi^{-1}(gK) \in M$ and hence that $a = t$ and $\z(\tilde \g) \subset \t$.
\end{pf}
From Lemma \ref{center} and Remark \ref{remark2.3}, we immediately obtain the following corollary.\par
\begin{cor} \label{asymptoticdeformation}  An  homogeneous bundle  $M= G^\C \times_{P, \tau}  F$ admits a K\"ahler-Ricci soliton
if and only if there is a K\"ahler-Ricci soliton $(\omega, X)$ on $M$, where
$\omega$ is  a $G \times T^m$-invariant   and  $X = J \hat Y$ for some $Y\in \t = Lie(T^m)$.
\end{cor}
\bigskip
\bigskip
\section{The Tian-Zhu  invariants of a  homogeneous toric bundles}
\bigskip
In \cite{TZ1},  G. Tian and X.-H. Zhu proved that, on a given compact complex manifold $(N, J)$, a vector field $X$ is the associated vector field  of  a K\"ahler-Ricci soliton  $(\omega, X)$
only if a certain   holomorphic invariant homomorphism
$$\Fut_X: \aut(N,J) \to \R$$
vanishes identically. Such  homomorphism  $\Fut_X$ is an analogue  of the classical Futaki invariant $\Fut:\aut(N,J) \to \R$ of $(N,J)$ (\cite{Fu}) and one has $\Fut_X = \Fut$ when $X = 0$. In the following, we will  call  such homomorphism   the {\it Tian-Zhu invariant
associated with $X$\/}. \par
\par
In \cite{TZ1} the following important property has been proved.\par
\medskip
\begin{theo} (\cite{TZ1}, Prop. 2.1) \label{TZtheorem1} Let $(N,J)$ be a compact complex K\"ahler manifold with $c_1(N) > 0$. For any maximal compact
subgroup $\tilde G \subset \Aut^o(N,J)$,  there exists exactly one element $Y \in \z(\tilde \g)$ (possibly equal to $0$) so that
$\Fut_{J \hat Y}( \cdot )$  vanishes identically.
\end{theo}
\bigskip
Let us now consider the toric bundle $M$. By Lemma \ref{center} and Theorem \ref{TZtheorem1},
if we consider a maximal compact subgroup $\tilde G \subset \Aut(M, J)$ that
contains $G \times T^m$, there exists exactly one $Y \in \t$ so that the Tian-Zhu invariant  $\Fut_X(\cdot)$ with
$X = J \hat Y$ vanishes.\par
\smallskip
We need now to determine the explicit expression for $\Fut_X$ when $X = J \hat Y$, for some $Y \in \t$.  Recall that, since $(Z'_1, \dots, Z'_\co)$ is a
basis for $\t$, any  vector field of this kind is of the form
$$X^{(\lambda)} = \sum_{\ell=1}^\co\lambda^\ell J \hat Z'_\ell \eqno(4.1)$$
 for some
suitable $\lambda = (\lambda^1, \dots,  \lambda^\co) \in \R^\co$.\par
\medskip
\begin{lem} \label{functiontheta} Assume that  $M$ is Fano and let $\omega$ be
a $G \times T^\co$-invariant K\"ahler form on $M$, with algebraic representative  so that
$Z_\omega|_F = \sum_i f_i Z_i + Z_V$ for some smooth functions  $f_i: F \to \R$. For any vector field
$X$ on $M$ such that $\L_{J X} \omega = 0$,  then there exists a  unique  smooth real valued function $\theta^{(X)}$ such that
$$\left\{ \begin{matrix}
\L_{X}ÿ\omega = \frac{1}{4 \pi} d d^c \theta^{(X)} = \frac{i}{2 \pi} \partial \bar \partial \theta^{(X)}\\
\phantom{aaaa}\\
\int_M e^{\theta^{(X)}} \omega^{n+\co} = \int_M \omega^{n+\co}\ .
\end{matrix}
 \right.
\eqno(4.2)$$
If $X = X^{(\lambda)}$   is a vector field of the form
(4.1), then the corresponding function $\theta^{(\lambda)}$ is $G \times T^\co$-invariant
and  the restriction of $\theta^{(\lambda)}|_F$ is
$$\left.\theta^{(\lambda)}\right|_F = - 4 \pi \sum_j\lambda^j f_j + C^{(\lambda)}\eqno(4.3)$$
where $C^{(\lambda)}$ is the real number
$$C^{(\lambda)} = \log\left(\frac{\int_M \omega^{n+\co}}{\int_M e^{-4 \pi \lambda^i f_i} \omega^{n + \co}}\right)\ .\eqno(4.4)$$
The constant $C^{(\lambda)}$ is the same for all  cohomologous $G \times T^\co$-invariant  K\"ahler forms.
\end{lem}
\begin{pf} Since $c_1(M) > 0$  and  $\L_{JX} \omega = d \imath_{JX}ÿ\omega = 0$,
by Bochner's theorem $b_1(M) = 0$ and  there exists a unique function $\theta_X$ so that
$$ \frac{1}{4 \pi} d^c \theta^{(X)} =  (\imath_{X}ÿ\omega) \circ J =  - \imath_{J X}ÿ\omega $$
and $(4.2)_2$ is satisfied. From uniqueness and the fact that $\omega$ and $X^{(\lambda)}$ are both $G \times T^\co$-invariant, it follows that  the function $ \theta^{(\lambda)}$ associated with $X = X^{(\lambda)}$ is $G \times T^\co$-invariant.  Moreover, one can  check   the $G \times T^\co$-invariant function defined by (4.3)   is the required function because
it  satisfies $(4.2)_2$ and
$$\left.d d^c \theta^{(\lambda)}(\hat Z_i, J \hat Z_j)\right|_F = \left.\L_{X^{(\lambda)} }\omega(\hat Z_i, J \hat Z_j)\right|_F \, \qquad 1 \leq i,j \leq \co\ .$$
 Finally, by (2.6) if $\omega' $ and $\omega$ are cohomologous,  the algebraic representative of $\omega'$ is given by
 $Z_{\omega'}|_F = \sum_j f'_j Z_j + Z_V$ with $f'_j = f_j - J \hat Z_i\left(\frac{1}{4 \pi} \phi\right)$ for
 some smooth $G \times T^\co$-invariant function $\phi: M \to \R$. By (4.3), the function $\theta'{}^{(\lambda)}$ relative to  $\omega'$ is
 $$\theta'{}^{(\lambda)} = - 4 \pi \sum_j\lambda^j f_j  + \lambda^j J \hat Z_j\left(\phi\right) + C'{}^{(\lambda)} =
 \theta^{(\lambda)}ÿ+ X^{(\lambda)}(\phi) + \left(C'{}^{(\lambda)} - C^{(\lambda)} \right)\ ,$$
 where we denoted by $C'{}^{(\lambda)} $ the constant  (4.4) determined by $\omega'$ in place of $\omega$.
 On the other hand,  Lemma 2.1 of  \cite{TZ1} shows  that $\theta'{}^{(\lambda)} =  \theta^{(\lambda)}ÿ+ X^{(\lambda)}(\phi)$ and  hence that  $C'{}^{(\lambda)} =  C^{(\lambda)}$.
 \end{pf}
\bigskip
\begin{prop} \label{TZinvariant}  Assume  that  $M$ is Fano. For any $\lambda \in \R^\co$, we have that
 $\Fut_{X^{(\lambda)}}(J\hat Y) = 0$ for all  $Y \in \t$ if and only if  all integrals
$$ \int_{\Delta_F}  x_k e^{ \lambda^a  x_a }
\prod_{\alpha \in R^+_\m}
\left( - \frac{a^j_{\alpha} x^j}{4 \pi}  + b_\alpha\right)
dx^1 \wedge \dots \wedge dx^\co\ ,\ \ 1 \leq k \leq \co\ ,\eqno(4.5) $$
vanish.
\end{prop}
\begin{pf}  Let $\omega$ be a $G \times T^\co$-invariant K\"ahler form  in  $c_1(M)$, $\rho$ the Ricci form of $\omega$  and  $h_\omega$
a smooth function on $M$ such that
$$\rho - \omega = \frac{1}{4 \pi} d d^c h_\omega = \frac{i}{2 \pi} \partial \bar \partial h_\omega\ .\eqno(4.6)$$
It follows that  the algebraic representatives of $\rho$, $\omega$ and $\frac{1}{4 \pi} d d^c h_\omega$   are so that  $Z_\rho - Z_\omega = Z_{\frac{1}{4 \pi} d d^c h_\omega}$. From
(2.6) - (2.8), we have that
$$J \hat Z_j(h_\omega)  =  - J \hat Z_j(\log h)  + 4 \pi  f_j \eqno(4.7)$$
where $h: F \to \R$ is defined in (2.9).  According to the definition  given in \cite{TZ1},
$$\Fut_{X^{(\lambda)}}(Y) \= \int_M Y\left( h_\omega - \theta^{(\lambda)}\right) e^{\theta^{(\lambda)}} \omega^{n+\co} \qquad
\text{for any}\ Y \in \aut(M, J)\ ,\eqno(4.8)$$
where   $\theta^{(\lambda)}: M \to\C$ is the unique smooth function that satisfies the condition (4.2).\par
 \medskip
 Let us now compute $\Fut_{X^{(\lambda)}}(J \hat Y)$ when $Y =  Z_k'$.  First of all, let us  fix a point $p_o \in F_\reg$ and consider the diffeomorphism
 $$\xi: F_\reg \to \R^\co \times T^\co \simeq (\C_*)^{\co}\ ,$$
 $$\xi\left(\exp\left(\sum_{j= 1}^\co \left(t^j + i s^j\right)i Z'_j\right) \cdot p_o\right) = \left(\frac{2 \pi}{\nu_1}  t^1 e^{i   \frac{2 \pi}{\nu_1} s^1}, \dots,
\frac{2 \pi}{\nu_\co}  t^\co e^{i \frac{2 \pi}{\nu_\co}
 s^\co}
 \right)\ .$$
If we identify $F_\reg$ with $ (\C_*)^{\co}$ by means of $\xi$,  the pairs $(\frac{2 \pi}{\nu_i} t^i,  \frac{2 \pi}{\nu_i} s^i)$'s  are
polar coordinates for  the factors  $\C_*$ of $(\C_*)^\co$ and we may consider the $\co$-tuple
$(\frac{2 \pi}{\nu_1}(t^1 + i s^1), \dots, \frac{2 \pi}{\nu_\co}(t^\co + i s^\co))$ as a system of complex
 coordinates on $F_\reg \simeq (\C_*)^\co$
 such that
 $$\frac{\partial}{\partial t^i} =   J \hat Z_i\ ,\qquad  \frac{\partial}{\partial s^i} = - \hat Z_i\ .$$
 Now, set $\Omega_F = dt^1 \wedge \dots \wedge d t^m \wedge ds^1 \wedge \dots \wedge d s^\co$. \par
 \begin{lem} \label{Fubini}
 There is a suitable constant $C$ such that for any $G$-invariant function
 $\phi\in \CC^\infty(M)^G$ we have
 $$\int_M \phi\cdot \omega^{n+m} = C \cdot \operatorname{Vol}_{\omega_V^n}(V) \cdot \int_{F_\reg} \phi\cdot  h \cdot \Omega_F\ ,\eqno(4.9)$$
 where $h \in \CC^\infty(F)^{T^\co}$ is the  $T^m$-invariant function defined in (2.9).  \end{lem}
 \begin{pfnqed} Let $\CU_\m$ be a open set in $\m$ containing $0$ such that the map
 $\psi: \CU_\m\to \exp(\CU_\m)\cdot (eP) \= \CU_V$ is a diffeomorphism onto its image and the mapping
 $\lambda:F\times \CU_V\to \pi^{-1}(\CU_V) \= \CU$ given by $\lambda(f,\exp(X)\cdot (eP)) = \exp(X)\cdot f$ is a
 bundle isomorphism.
We then select $g_j\in G$, $j=1,\ldots,N$, so that $V=\bigcup_{j=1}^N (g_j\cdot \CU_V)$ and put $A_0=\emptyset$ and
$A_j \= g_j\CU$ for $j=1,\ldots,N$. Hence
  $$\int_M  \phi \cdot \omega^{n+m} = \sum_{j=1}^N \int_{A_j\setminus \left(\bigcup_{i=0,}^{j-1}A_i\right)}  \phi
  \cdot \omega^{n+m}\ =$$
  $$ = \sum_{j=1}^N \int_{\CU\setminus \left(\bigcup_{i=0,}^{j-1}g_j^{-1}A_i\right)}  \phi
  \cdot \omega^{n+m}\ .\eqno (4.10)$$
In $\CU\cong F\times \CU_V$ we may restrict to the submanifold $F_{\reg}\times \CU_V$ and we may define the function $\tilde h
\in C^\infty(F_\reg\times \CU_V)$ by means of the following
 $$\omega^{n+m} = \tilde h \cdot \Omega_F\wedge \omega_V^n\ . $$
The function $\tilde h$ can be easily determined  by evaluating the forms $\omega^{n+m}$ and
$\Omega_F\wedge \omega_V^n$ on the frames $\{\hat Z_i, J\hat Z_i, \hat F_j, \hat G_j\}$ at the points of $F_\reg$; A
direct computation shows that $\tilde h = C \cdot h$, for some suitable constant $C$. \par
Using the $G$-invariance and Fubini's theorem,  (4.10) reads
 $$\int_M  \phi \cdot \omega^{n+m}  = \sum_{j = 1}^N \int_{F\times \left(\CU_V\setminus
 \left(\bigcup_{i=1}^{j-1}g_j^{-1}g_i\CU_V\right)\right)}  \phi \cdot \omega^{n+m} =$$
$$= \sum_{j = 1}^N \int_{\left(\CU_V\setminus
 \left(\bigcup_{i=1}^{j-1}g_j^{-1}g_i\CU_V\right)\right)}\left(\int_{F_{\reg}} \phi \cdot \tilde h \cdot \Omega_F\wedge
\omega_V^n  \right) =$$
$$ =  C \left(\sum_{j = 1}^N \int_{\CU_V\setminus
 \left(\bigcup_{i=1}^{j-1}g_j^{-1}g_i\CU_V\right)} \omega_V^n\right) \cdot \int_{F_{\reg}} \phi\cdot h \cdot \Omega_F =
C \cdot \operatorname{Vol}_{\omega_V^n}(V) \cdot \int_{F_\reg} \phi\cdot  h \cdot \Omega_F\ .\eqno\blacksquare$$
 \end{pfnqed}
 By Lemma \ref{Fubini}, (4.8), (4.2) and (4.7), it follows that  $\Fut_{X^{(\lambda)}}(J\hat Y) = 0$ for all  $Y \in \t$ if and only if  the integrals
  $$ g^{(\lambda)}_k = \int_{F_\reg}   J \hat Z_k\left( h_\omega - \theta^{(\lambda)}\right) e^{\theta^{(\lambda)}} \cdot h \cdot
 \Omega_F = $$
 $$ = \int_{F_\reg} \left(
- J \hat Z_k (\log  h)   +  4 \pi f_k  +  J \hat Z_k(   4 \pi  \lambda^j  f_j)\right) \cdot
e^{- 4 \pi \lambda^j f_j + C^{(\lambda)}} \cdot h \cdot \Omega_F\eqno(4.11)$$
are equal to $0$ for all $k = 1, \dots, \co$.\par
On the other hand,  if we identify $F_\reg$  with $(\C_*)^\co$ by means of  the map $\xi$ described above, we have that
$$\int_{F_\reg} \left(
J \hat Z_k (\log h)     -  J \hat Z_k(  4 \pi  \lambda^j  f_j)\right) \cdot
e^{- 4 \pi \lambda^j f_j} \cdot h \cdot \Omega_F = $$
$$ =
\int_{\R^\co}
\frac{\partial}{\partial t^k}  \left( h  \cdot
e^{- 4 \pi  \lambda^j f_j} \right) dt^1 \wedge \dots \wedge d t^\co \wedge ds^1 \wedge \dots \wedge
ds^\co = $$
$$ =  (-1)^{k}\int_{\R^{\co - 1}}  \left(\int_{-\infty}^\infty \frac{\partial}{\partial t^k}  \left( h  \cdot
e^{- 4 \pi  \lambda^j f_j} \right) dt^k\right)  dt^1 \wedge \dots \underset {k}{\widehat{\phantom{A}}}\dots \wedge d t^\co \wedge ds^1 \wedge \dots \wedge
ds^\co = $$
$$ = (-1)^{k}\!\!\! \int_{\R^{\co - 1}}  \left(\lim_{\smallmatrix a \to + \infty\\
b \to - \infty \endsmallmatrix}  \left.\left( h  \cdot
e^{- 4 \pi  \lambda^j f_j} \right) \right|_{t^k = a}^{t^k = b} \right)  dt^1 \wedge \dots \ \underset {k}{\widehat{\phantom{A}}}\dots \wedge d t^\co \wedge ds^1 \wedge \dots \wedge
ds^\co  = $$
$$ = 0\eqno(4.12)$$
where  the last equality is obtained from the definition of $h$,  the fact that the functions $f_j: F_\reg \to \R$ are bounded (see (2.12))   and the property that
$$\det(- f_{i,j}): F_\reg \to \R\ ,\qquad \left.\det(- f_{i,j})\right|_p = \frac{1}{(2 \pi)^\co} \det(g_p(\hat Z_i, \hat Z_j)) $$
goes  to $0$ when $p$ tends to a point of $F \setminus  F_\reg$, since
$$F \setminus  F_\reg = \{\ q \in F\ \ \text{such that}\ \left.\hat Z_j\right|_q = 0\ \text{for some}\ j = 1, \dots, \co\ \}\ .$$
From (4.12), it follows that the integrals  $g^{(\lambda)}_k $ are equal to
$$g^{(\lambda)}_k=  C \int_{F_\reg}
f_k
e^{- 4 \pi  \lambda^j f_j}  \det(- f_{\ell,s}) \cdot \prod_{\alpha \in R^+_\m}
(a_{ \alpha}^i  f_j + b_\alpha) \cdot \Omega_F \eqno(4.13)$$
for some constant $C$.
Using the change of variables $(t^i, s^i) \mapsto (x^i = - 4 \pi f^i(t^j,s^k), s^i)$ and the the fact that the integrand
is independent of the
coordinates $s^i$,
one can check that (4.13) is equal (up to a multiplicative constant) to the integral (4.5) over  the image of the moment
map $\mu = (- 4 \pi   f_1, \dots, - 4 \pi   f_\co)$, i.e. the canonical polytope  $\Delta_F \subset \R^\co$. \end{pf}
\bigskip
\bigskip
\section{The reduction of the solitonic K\"ahler equation on $M$ to an equation on the toric manifold $F$ and the proof of Theorem \ref{maintheorem}} \label{reduction}
\bigskip
By Thm. 1.1 of \cite{PS2}, we know  that    $M = G^\C \times_{P, \tau} F$ is Fano only if  also $F$ is Fano. In the proof of that theorem, we have shown that the correspondence $\omega \mapsto \omega|_{TF}$
between 2-forms on $M$ and on $F$ maps
any $G$-invariant K\"ahler form  in  $c_1(M)$ into a   $T^\co$-invariant K\"ahler form  on $F$, which belongs  to $c_1(F)$. The following lemma shows that such correspondence  is actually bijective.\par
\medskip
\begin{lem} \label{reductionlemma}ÿLet $c_1(M) > 0$ and denote by $c_1(M)^{G\times T^\co}$ and $c_1(F)^{T^\co}$ the sets of $G\times T^\co$-invariant 2-forms in $c_1(M)$ and of
$T^\co$-invariant 2-forms in $c_1(F)$, respectively. Then there exists a map
$$E:\ \  c_1(F)^{T^\co} \longrightarrow c_1(M)^{G\times T^\co}  \eqno(5.1)$$
which is inverse to the map
$$R:\ \ \omega \in c_1(M)^{G\times T^\co} \quad  \longrightarrow \quad \omega|_{TF} \in c_1(F)^{T^\co}\ .$$
Moreover, $E( \omega) $ is K\"ahler if and only if $\omega $ is K\"ahler.
\end{lem}
\begin{pf} Let us fix a K\"ahler form $\omega_o \in c_1(F)^{T^\co}$.   Any
 $T^\co$-invariant  $\omega \in c_1(F)$ is of the form $\omega = \omega_o + d d^c \phi_\omega$ for
 some $T^\co$-invariant function $\phi_\omega$, which is unique up to a constant.
We denote by $\mu_\omega$ the map
 $$\mu_\omega  : F \to  \t^*\ ,\qquad \mu_\omega|_p(X) \= \mu_{\omega_o}|_p(X) - d^c \phi_\omega(\hat X_p)
 \ \text{for any}\ \ X \in \t .\eqno(5.2)$$
 where $\mu_{\omega_o}$ is the metrically normalized moment
map associated to the K\"ahler form $\omega_o$ (for the definition, see \S 2).  One can check that
$\mu_\omega$ is the metrically normalized moment map relative to $\omega$, whenever $\omega$
is non-degenerate. \par
\smallskip
Now, for any $\omega \in c_1(F)^{T^\co}$  we define
 $E(\omega)$ as the unique $G\times T^\co$-invariant 2-form on $M$,  whose restriction on  $TM|_F$
is as follows: for any $p\in F$, $X, Y \in T_pF$ and $A, B \in \m$
$$E(\omega)_p(X, Y) = \omega_p(X, Y)\ ,\qquad E(\omega)_p(X, \hat A) = 0\ ,$$
$$E(\omega)_p(\hat A, \hat B) = - \mu_\omega(p)(\tau([A,B]_\k)) + (\pi^*\omega_V)_p(\hat A, \hat B)\ ,\eqno(5.3)$$
where $\omega_V$ is the $G$-invariant K\"ahler-Einstein form on $V$ with Einstein constant $c = 1$ and where  we denoted by  ``$[A,B]_\k$"  the component of $[A,B]$  along $\k$ w.r.t. the decomposition
$\g = \k \oplus \m$. Going through  the arguments  after formula (5.1) of  \cite{PS2}, one can check that
$E(\omega)$  is closed and $J$-invariant. \par
It is also direct to see that the algebraic representative $Z_{E(\omega)}$ of $E(\omega)$ is so that
$$\left.Z_{E(\omega)}\right|_F = \sum_j \left(f_{oj} - J \hat Z_j(\phi_\omega)\right) Z_j + Z_V\eqno(5.4)$$
  where the $f_{oj}: F \to \R$ are (up to the factor $- 4 \pi$) equal  to the components of $\mu_{\omega_o}: F \to \t^* \simeq \R^\co$ under the
  identification (2.12). It follows that the algebraic representative of $E(\omega_2) - E(\omega_1)$
  is the same of $d d^c (\phi_{\omega_2} - \phi_{\omega_1})$, meaning that the image of $E$ is in a
  single cohomology class.
  Moreover, by looking at the algebraic representatives, one can see that the 2-form $E (\rho)$, where $\rho$
is the Ricci form of $\omega$,    coincides with
 the 2-form $\rho_o $ defined in formula (5.7) of \cite{PS2}. By the proof of Thm. 1.1 in \cite{PS2},
 we know that $\rho_o \in c_1(M)$ and hence $E(\omega) \in c_1(M)$ for any $\omega$.\par
 \smallskip
From (5.4) and the remarks at the end of \S 2,  for any $\omega \in   c_1(M)^{G\times T^\co}$
the algebraic representatives of $\omega$ and of $E(R(\omega))$ coincide and hence
$\omega = E(R(\omega))$. This implies that $E$ is inverse to $R$ since,  by construction,  we also have that $R(E(\omega)) = \omega$ for any $\omega \in c_1(F)^{T^\co}$. \par
\smallskip
The last claim follows from the fact that, for any K\"ahler form $ \omega \in c_1(F)^{T^\co}$, the 2-form
$E(\omega)$ is positive because it is $G$-invariant and
 its restriction at  $TM|_F$ is positive. This is true because,  if we denote by $- 4 \pi f_j$ the
components of the metrically normalized moment map $\mu_\omega$ of $\omega$, from (5.4)
we have that
  for any $\alpha_j \in R^+_\m$
$$\tilde \omega_p(\hat F_{\alpha_j},J \hat F_{\alpha_j}) = \B(\sum_k f_k Z_k + Z_V, [F_{\alpha_j},
G_{\alpha_i}]) = i \alpha_j(\sum_k f_k Z_k + Z_V) > 0\ ,\eqno(5.5)$$
 which coincides with condition (1.1)   of Thm. 1.1 of \cite{PS2} in a different notation.
\end{pf}
 \bigskip
We now want to determine the differential equations that characterize the
 $T^\co$-invariant K\"ahler forms  in $c_1(F)$  corresponding to invariant solitonic K\"ahler  forms of $M$. We  recall that, by Corollary \ref{asymptoticdeformation}, there exists a solitonic
 K\"ahler form  on $M$ if and only if there exists a K\"ahler-Ricci soliton $(\tilde \omega, X)$ where
 $\tilde \omega$ is $G \times T^\co$-invariant and  $X$ is of the form $X = X^{(\lambda)} = \sum_k\lambda^k J \hat Z'_k$ for some $\lambda \in \R^\co$.\par
 \medskip
 To simplify the notation, in the following,
for any $T^\co$-invariant 2-form $\omega \in c_1(F)$, we will denote by $\tilde \omega = E(\omega)$ the
 corresponding  2-form in $c_1(M)^{G \times T^\co}$. Similarly, any 2-form   in $c_1(M)^{G \times T^\co}$ will be denoted with a symbol of the form $\tilde \omega$ and the corresponding 2-form in $c_1(F)$
 will be indicated by $\omega = R(\tilde \omega)$.\par
 \medskip
 For any  K\"ahler form $\tilde \omega \in  c_1(M)^{G \times T^\co}$, let us indicate by $(- 4 \pi) \cdot f_{\omega,j}: F \to \R$ the components of the
  metrically normalized moment map $\mu_\omega: F \to \Delta_F \subset \R^\co \simeq \t^*$, relative to $\omega$, as given  in (2.11) and (2.12). Let also $\tilde \rho$ be the Ricci form of $\tilde \omega$ and
 $\phi_{\tilde \omega, \tilde \rho}$  the unique potential on $M$ so that
 $$\tilde \rho = \tilde \omega + \frac{1}{4 \pi}ÿd d^c \phi_{\tilde\omega, \tilde\rho} \ , \qquad \quad
 \int_M e^{\phi_{\tilde\omega, \tilde\rho}} \tilde \omega^{n + \co} =  \int_M \tilde \omega^{n + \co}.\eqno(5.6)$$
 Notice that $\phi_{\tilde \omega, \tilde \rho}$ is $G \times T^\co$-invariant and hence it is uniquely determined by its restriction to $F$.\par
 From  Lemma \ref{functiontheta}, the pair $(\tilde \omega, X^{(\lambda)})$ is
 a K\"ahler-Ricci soliton if and only  if
 $$\phi_{\tilde \omega, \tilde \rho} + 4 \pi \sum_\lambda^ij f_{\omega,j} -  C^{(\lambda)} = 0\ ,\eqno(5.7)$$
where $C^{(\lambda)}$ is defined by (4.4) and it is independent of $\omega$.\par
\medskip
We now fix a $T^\co$-invariant K\"ahler form $\omega_o \in c_1(F)$. For any other $T^\co$-invariant
K\"ahler form $\omega$, we denote by  $\psi_\omega$ the unique potential on $F$ so that
$$\omega = \omega_o + \frac{1}{4 \pi}ÿd d^c \psi_{\omega}\ ,\qquad \quad
 \int_F e^{\psi_{\omega}} \omega^{\co}_o =  \int_F \omega^{\co}_o\eqno(5.8)$$
and we  want to determine the equation in the unknown function $\psi_\omega$ determined by the condition (5.7).\par
\bigskip
Let $\tilde \omega_o = E(\omega_o)$ be the $G \times T^\co$-invariant K\"ahler form given by Lemma
\ref{reductionlemma} and   let
  $Z_{\tilde \omega_o} = \sum_i f_{oi} Z_i + Z_V$ be the restriction to $F$ of the algebraic representative
  of $\tilde \omega_o$.   We consider also
a system of complex coordinates  $(t^1 + i s^1, \dots, t^\co + i s^\co)$ on $F_\reg \simeq (\C_*)^{\co}$
 as in the proof of Proposition \ref{TZinvariant},  such that
 $$\frac{\partial}{\partial t^i} = J \hat Z_i\ ,\qquad  \frac{\partial}{\partial s^i} = - \hat Z_i\ .\eqno(5.9)$$
Since $F_\reg \simeq  \R^{ \co} \times T^\co$, the maps $f_{oi}$ are $T^\co$-invariant
  and, by  (2.10),  $\frac{\partial f_{oi}}{\partial t^j} =  \frac{\partial f_{oj}}{\partial t^i}  $ for all $i,j$'s, then there exists a $T^\co$-invariant smooth function  $u_o: F_\reg  \to \R$ so that
 $$\left.f_{oi}\right|_{F_\reg} = - \frac{1}{4 \pi} \frac{\partial u_o}{\partial t^i}\ .\eqno(5.10)$$
 The function $u_o$ is uniquely determined up to an additive constant. It can also be checked that
 $\omega_o|_{F_\reg} = \frac{1}{4 \pi} d d^c u_o$.\par
We claim that there exists some suitable constant $C_\omega$ so that
$\phi_{\tilde\omega, \tilde\rho}|_{F_\reg} = \Psi + C_\omega$ where
$$\Psi \= - \log \left|  \det\left(\frac{\partial^2  (u_o + \psi_\omega)}{\partial t^i \partial t^j}ÿ\right)\cdot \prod_{\alpha \in R^+_\m} \left( - \frac{a^i_\alpha}{4 \pi}
 \frac{\partial (u_o + \psi_\omega)}{\partial t^i} + b_\alpha\right)\right| -  (u_o  +  \psi_\omega)\ .\eqno(5.11)$$
In order to check this, notice that  by (2.6), the restriction to $F$ of the algebraic representative of $\tilde \omega = \tilde \omega_o + \frac{1}{4 \pi} d d^c \psi_\omega$ is
 $$Z_{\tilde \omega} = \sum_j \left(f_{oj} - \frac{1}{4 \pi}ÿJ Z_j(\psi_\omega)\right)Z_j  + Z_V =
 - \frac{1}{4 \pi} \sum_j \frac
 {\partial (u_o + \psi_\omega)}{\partial t^j} Z_j + Z_V\ .$$
 Then, from (2.8), it follows that the algebraic representative of $\frac{1}{4 \pi} d d^c \Psi$
  coincides with the
 algebraic representative of $\tilde \rho - \tilde \omega = \frac{1}{4 \pi} d d^c \phi_{\tilde \omega, \tilde \rho}$ and hence that $ \Psi -  \phi_{\tilde \omega, \tilde \rho}$ is a constant. The value of  $C_\omega$ is uniquely determined  by
 the normalizing condition $(5.6)_2$. \par
 \smallskip
 From (5.7) together with the expression of $\phi_{\tilde\omega, \tilde\rho}$ given by (5.11) and  setting  $\varphi = \psi_\omega + C_\omega$, we obtain  the following proposition, which reduces the solitonic  equation on $M$ to a Monge-Ampere equation on the toric manifold  $F$.\par
\medskip
\begin{prop} \label{solitonicequation}ÿLet $\omega_o$ be a fixed $T^\co$-invariant K\"ahler form in $c_1(F)$ and $u_o: F_\reg \to \R$
be a fixed smooth function so that (5.10) holds.  Let also $X^{(\lambda)}$ be the  vector field of the form (4.1) such that $\Fut_{X^{(\lambda)}}(\cdot) = 0$ and let $C^{(\lambda)}$ be the constant defined by (4.4). \par
A $T^\co$-invariant K\"ahler metric $\omega \in c_1(F)$ is so that the corresponding $G \times T^\co$-invariant K\"ahler form $\tilde \omega \in c_1(M)$ is a soliton, with associated vector field $X^{(\lambda)}$, if and only if
$$\omega = \omega_o + \frac{1}{4 \pi}ÿd d^c \varphi\ ,$$
 where $\varphi: F \to \R$ is a smooth $T^\co$-invariant function that satisfies
$$ \det\left(\frac{\partial^2  (u_o + \varphi)}{\partial t^i \partial t^j}ÿ\right) =
\frac{1}{\prod_{\alpha \in R^+_\m} \left( - \frac{a^i_\alpha}{4 \pi}
 \frac{\partial (u_o + \varphi)}{\partial t^i} + b_\alpha\right) }e^{ - C^{(\lambda)} -   X^{(\lambda)} \left(u_o + \varphi\right) -  u_o - \varphi} \eqno(5.12)$$
 at all points of $F_\reg$.
\end{prop}
\bigskip
We recall that,  by Theorem \ref{TZtheorem1} and Lemma \ref{center}, there exists a unique $\lambda \in \R^\co$
 so that  $\Fut_{X^{(\lambda)}}(\cdot)$ vanishes identically.  By Proposition \ref{TZinvariant}, for such $\lambda$ all integrals (4.5) are equal to $0$. \par
 \medskip
 In \cite{WZ},  X.-J. Wang and X. Zhu determined a Monge-Ampere equation that
  characterizes the $T^\co$-invariant K\"ahler-Ricci
 solitons on   $F$ and proved the solvability of such equation.  Notice that the case considered
 by Wang and Zhu can be interpreted as a homogeneous toric bundle with  basis given by a single point.
 And in fact  the
 Monge-Ampere of
  Wang and Zhu  can be obtained from
(5.12)  by  setting the factor
  $${\mathcal A} =  \frac{1}{\prod_{\alpha \in R^+_\m} \left( - \frac{a^i_\alpha}{4 \pi}
 \frac{\partial (u_o + \varphi)}{\partial t^i} + b_\alpha\right) }\eqno(5.13)$$
 equal to $1$.  \par
We  claim  that  the  arguments used  in \cite{WZ} for proving   the
solvability of
 (5.12)  when  ${\mathcal A} = 1$  remain valid also when ${\mathcal A} \neq 1$  and hence
 that (5.12) is always solvable. \par
 \bigskip
 In fact, as in \cite{WZ},  the solvability of (5.12) can be obtained by  the continuity method, namely by
 considering  the family of equations
 $$ \det\left(\frac{\partial^2  (u_o + \varphi)}{\partial t^i \partial t^j}ÿ\right) =
\frac{1}{\prod_{\alpha \in R^+_\m} \left( - \frac{a^i_\alpha}{4 \pi}
 \frac{\partial (u_o + \varphi)}{\partial t^i} + b_\alpha\right) }e^{ - C^{(\lambda)} -   X^{(\lambda)} \left(u_o + \varphi\right) -  u_o - t \varphi}\ , \eqno(5.14)$$
parameterized by the real numbers  $t \in [0,1]$. By the same arguments used for the proof of Proposition \ref{solitonicequation}, one can see that   (5.14) is the Monge-Ampere equation characterizing the $T^\co$-invariant 2-forms $\omega \in c_1(F)$, whose associated $G \times T^\co$-invariant 2-form $\tilde \omega$ satisfy
 $$\tilde \rho - \tilde \omega_o - t (\tilde \omega - \tilde \omega_o) = \L_{X^{(\lambda)}} \tilde \omega\ .\eqno(5.15)$$
 By the results of \cite{Zh} and \cite{TZ}, (5.15)   is solvable for any $t$ in an open subinterval $[0, \epsilon[ \subset [0,1]$ and hence the same is true for the equation (5.14). Moreover,
 for any $t \in [0,1]$, if  $\varphi$ is a solution of (5.14) that corresponds  to a K\"ahler metric, then its
$G \times T^\co$-invariant extension is the potential w.r.t. $\tilde \omega_o$ of a  K\"ahler form $\tilde \omega$ that satisfies (5.15) and hence  in $c_1(M)$. It follows that the functions
 $$- 4 \pi f_i =  \frac{\partial (u_o + \varphi)}{\partial t^i}$$
 are the components of the metrically normalized moment map of $\omega$ and take value in the polytope $\Delta_F$, which is a bounded convex domain in $\R^\co$ and it is independent of
 $t$.  In particular, the algebraic representative of $\tilde \omega$ is of the form (2.7),  the
 function $\mu = (- 4 \pi f_1, \dots, - 4 \pi f_\co): F \to \R^\co \simeq t^*$ is a  metrically normalized moment map on $F$ and the image $\mu(F) = \Delta_F$ is independent of the value of  $t$.  Since $\mathcal A$ coincides with
 the function $\prod_{\alpha \in R^+_\m} \left( - \B(\sum_{i = 1}^\co f_i Z_i + Z_V, H_\alpha\right) $, by Thm 1.1
 of \cite{PS2}, the values of $\mathcal A$ are always positive and are bounded above and below
 by two constants that are independent of $t$.\par
 From this facts and  the generalizations of  \cite{TZ} of the a priori estimates of \cite{Ya1},
 the solvability of (5.14) for $t = 1$  is proved if one can give a  uniform upper and lower estimates
  for the solutions $\varphi$ of (5.14) for any $t \in [\epsilon, 1]$ for some $0 < \epsilon$.\par
 \medskip
 One can  check that  the proofs of Lemmata 3.2 -  3.5
 in \cite{WZ} remain
 valid also if ${\mathcal A} \neq 1$ provided that the following properties and differences on notation  are  taken into account:
 \begin{itemize}
 \item[i)] The canonical polytope $\Delta_F$ should be considered as equal to the  dual $\Omega^*$ of the polytope denoted by $\Omega$ in \cite{WZ} and the $m$-tuple $(\lambda^1, \dots, \lambda^\co)$
 should be considered as equal to the $m$-tuple of constants $(c_1, \dots, c_m)$
 considered in \cite{WZ}; Moreover, with  no loss of generality, one  should  assume that our function
  $u_o : F_\reg \to \R$ so that $\omega_o|_{TF} = \frac{1}{4 \pi}ÿd d^c u_o$ coincides with  the function denoted by ``$u^0$" in \cite{WZ};
 \item[ii)]ÿBy the previous remarks, if $M$ is Fano, there exist two positive real numbers $0 < K_1 < K_2$,
 independent on $t$, so that
 $K_1 \leq {\mathcal A} \leq K_2$ at any point of $F_\reg$ and for any solution of (5.12);
 \item[iii)] From (ii) and the fact that   all integrals (4.5) are equal to $0$,
 the polytope $\Omega^* = \Delta_F \subset \R^\co$ contains the origin also when $\mathcal A \neq 1$;
 \item[iv)] The
 equality $0 = \int_{\Omega^*} y_i e^{\sum_{\ell} c_\ell y_\ell} d y$ which appears at the end of the proof
 of Lemma 3.3 of \cite{WZ}  should be replaced by the equality
 $$0 = \int_{\Omega^*}ÿy_i e^{\sum_\ell c_\ell y_\ell} \prod_{\alpha \in R^+_\m}
\left( - \frac{\sum_j a^j_{\alpha} y^j}{4 \pi}  + b_\alpha\right)
dy^1 \wedge \dots \wedge dy^\co$$
which is true  by Proposition \ref{TZinvariant}; Under this replacement, all  remaining equalities considered in Lemma 3.3 of \cite{WZ} remain true also when ${\mathcal A} \neq 1$.
 \end{itemize}
 Those lemmata give the needed  estimates and  the solvability  of (5.12) is proved.\par

\bigskip
\bigskip
\font\smallsmc = cmcsc8
\font\smalltt = cmtt8
\font\smallit = cmti8
\hbox{\parindent=0pt\parskip=0pt
\vbox{\baselineskip 9.5 pt \hsize=3.1truein
\obeylines
{\smallsmc
Fabio Podest\`a
Dip. Matematica e Appl. per l'Architettura
Universit\`a di Firenze
Piazza Ghiberti 27
I-50100 Firenze
ITALY
}\medskip
{\smallit E-mail}\/: {\smalltt podesta@math.unifi.it
}
}
\hskip 0.0truecm
\vbox{\baselineskip 9.5 pt \hsize=3.7truein
\obeylines
{\smallsmc
Andrea Spiro
Dip. Matematica e Informatica
Universit\`a di Camerino
Via Madonna delle Carceri
I-62032 Camerino (Macerata)
ITALY
}\medskip
{\smallit E-mail}\/: {\smalltt andrea.spiro@unicam.it}
}
}

\end{document}